\documentclass[11pt]{amsart}

\usepackage{amscd}
\usepackage{amssymb}
\usepackage{boxedminipage}
\usepackage{enumerate}

\theoremstyle{plain} \newtheorem{theorem}{Theorem}[section]
\theoremstyle{plain} \newtheorem{lemma}[theorem]{Lemma}
\theoremstyle{plain} \newtheorem{proposition}[theorem]{Proposition}

\newtheorem{problem}[theorem]{Problem}

\newtheorem{conjectureP}{Conjecture}

\newcommand{\nr}{\refstepcounter{theorem}  
                   \noindent {\thetheorem .}}
\newcommand{\defi}{\medskip \noindent {\it Definition \nr} }
\newcommand{\defifin}{\medskip}
\newcommand{\eks}{\medskip \noindent {\it Example \nr} }
\newcommand{\eksfin}{\medskip}
\newcommand{\rem}{\medskip \noindent {\it Remark \nr} }
\newcommand{\remfin}{\medskip}

\newcommand{\coh}{{{\text{{\rm coh}}}}}

\newcommand{\modstab}[1]{{#1}-\underline{\text{mod}}}

\newcommand{\gT}{{\mathcal T}}
\newcommand{\gS}{{\mathcal S}}
\newcommand{\gF}{{\mathcal F}}
\newcommand{\gE}{{\mathcal E}}
\newcommand{\op}{{\mathcal O}}
\newcommand{\gFp}{{\mathcal F}^\prime}
\newcommand{\opw}{\op_{\psp{W}}}
\newcommand{\go}{\op}

\newcommand{\Hom}{\text{Hom}}

\newcommand{\gext}{\gE xt}

\newcommand{\im}{\text{im}\,}

\newcommand{\sus}{\subseteq}
\newcommand{\pil}{\rightarrow}
\newcommand{\rpil}{\leftarrow}

\newcommand{\vpil}{\leftarrow}
\newcommand{\inpil}{\hookrightarrow}

\newcommand{\projpil}{\dashrightarrow}
\newcommand{\dotpil}{\dashrightarrow}

\newcommand{\mto}[1]{\stackrel{#1}\longrightarrow}

\newcommand{\iso}{\cong}
\newcommand{\te}{\otimes}

\newcommand{\tL}{\tilde{L}}

\newcommand{\tV}{\tilde{V}}

\newcommand{\ome}{\omega_E}

\newcommand{\dl}{\Delta}
\newcommand{\cdel}{{C\Delta}}

\newcommand{\dlst}{\Delta^*}
\newcommand{\Sdl}{{\mathcal S}({\dl})}
\newcommand{\ord}{\succ}

\newcommand{\psp}[1]{{{\bf P}({#1})}}

\newcommand{\PP}{{\bf P}}
\newcommand{\psmi}{{\bf P}^{s-1}}
\newcommand{\pnmi}{{\bf P}^{n-1}}

\newcommand{\opnmi}{\go_{\pnmi}}

\newcommand{\hele}{{\bf Z}}

\def\ZZ{{\mathbb Z}}

\begin{document}
\title [BGG and simplicial complexes]
{(Bi-)Cohen-Macaulay simplicial complexes and their associated 
coherent sheaves}
\author { Gunnar Fl{\o}ystad \and Jon Eivind Vatne}
\address{ Matematisk institutt\\
          Johs. Brunsgt. 12 \\
          N-5008 Bergen \\
          Norway}   
        
\email{ gunnar@mi.uib.no \and  jonev@mi.uib.no }

\begin{abstract}
Via the BGG correspondence a simplicial complex $\dl$ on $[n]$ is
transformed into a complex of coherent sheaves on $\pnmi$.
We show that this complex reduces to a coherent sheaf $\gF$ exactly when
the Alexander dual $\dlst$ is Cohen-Macaulay.

We then determine when both $\dl$ and $\dlst$  are
Cohen-Macaulay.  This corresponds to $\gF$ being a locally
Cohen-Macaulay sheaf.

Lastly we conjecture for which range of invariants of such $\dl$'s it must be
a cone, and show the existence of such $\dl$'s which are not cones outside
of this range. 
\end{abstract}

\maketitle

\section*{Introduction}

To a simplicial complex $\Delta$ on the set
$[n] = \{1, \ldots, n\}$ is associated a monomial ideal $I_\Delta$ in the 
exterior algebra $E$ on a vector space of dimension $n$. Lately there has
been a renewed interest in the Bernstein-Gel'fand-Gel'fand (BGG) correspondence
which associates to a graded module $M$ over the exterior algebra $E$ 
a complex of coherent sheaves on the projective space $\PP^{n-1}$
(see \cite{EFS}, \cite{AAH99}). In this paper we study simplicial complexes in 
light of this correspondence.
Thus to each simplicial complex $\dl$ we get associated a complex of coherent
sheaves on $\PP^{n-1}$. Our basic result is that this complex reduces to
a single coherent sheaf $\gF$ if and only if the Alexander dual $\dlst$ 
is a Cohen-Macaulay simplicial complex.
So we in yet a new way establish the naturality of the concept of 
a simplicial complex being Cohen-Macaulay in addition to the well established
interpretations via the topological realization, and via commutative algebra
and Stanley-Reisner rings.

It also opens up the possibility to study simplicial complexes from the
point of view of algebraic geometry. A simple fact is that $\dl$ is a cone
if and only if the support of $\gF$ is contained in a hyperplane. 
Now the nicest coherent sheaves on projective space may be said to be 
vector bundles, or more generally those sheaves which when projected
down as far as possible, to a projective space of dimension equal to 
the dimension of the
support of the sheaf, become vector bundles. This is the class of locally
Cohen-Macaulay sheaves (of pure dimension). 
We show that the coherent sheaf $\gF$ is a locally Cohen-Macaulay sheaf 
iff both $\dl$ and $\dlst$ are Cohen-Macaulay simplicial complexes. We
call such $\dl$ bi-Cohen-Macaulay and try to describe this class as well
as possible.

In Section 1 we recall basic facts about the BGG-correspondence.
In Section 2 we apply this to simplicial complexes and show the basic
theorem, that we get a coherent sheaf $\gF$ via the BGG-correspondence iff
$\dlst$ is Cohen-Macaulay. We are also able to give a kind of geometric
interpretation of the $h$-vector of $\dlst$ in terms of the sheaf $\gF$.

In Section 3 we consider bi-Cohen-Macaulay simplicial complexes $\dl$.
Then the associated sheaf $\gF$ on $\PP^{n-1}$
when projected down to $\PP^{s-1}$, where $s-1$ is the dimension of the
support of $\gF$, becomes
one of the sheaves of differentials $\Omega^c_{\PP^{s-1}}$. This
gives quite restrictive conditions on the face vector of such $\Delta$. It is
parametrized by three parameters, namely $n$, $c$ and $s$.

When $c = 0$, $\Delta$ is just the empty simplex. When $c = 1$
a result of Fr\"oberg \cite{Fr2} enables a combinatorial 
description of such $\Delta$. If $\Delta$ has dimension $d-1$ equal to $1$ it
is a tree and in general $\dl$ is what is called a $(d-1)$-tree. 
When $c \geq 2$ a
combinatorial description seems less tractable. There is a classical example
of Reisner \cite{Re} of a triangulation $\Delta$ of the real projective plane
which is bi-Cohen-Macaulay if char $k \neq 2$;
but neither $\Delta$ nor $\Delta^*$ are Cohen-Macaulay if char $k = 2$.
In particular $\dl$ is not shellable.

Now suppose $\gF$ projects down to $\Omega^c_{\PP^{s-1}}$. A natural
question to ask is whether $\gF$ is degenerate or not (the support
contained in a hyperplane or not). This corresponds
to $\Delta$ being a cone or not.  In the last section, Section 5, we
conjecture that there exists a bi-Cohen-Macaulay  $\Delta$ which is not
a cone if and only if $s \leq n \leq (c+1)(s-c)$. 
We prove this conjecture when 
$c = 1$ and give examples to show the plausibility of this conjecture for 
any $c$. Also in the full range we prove the existence of bi-Cohen-Macaulay
$\Delta$ which are not cones.

As further motivation for the significance of bi-Cohen-Macaulay simplicial
complexes, we refer to the paper \cite{Fl2}. There a natural algebraically
defined family of simplicial complexes is defined which depends on the
parameters $n,d,c$, and an integer $a \geq 0$. It contains Alexander duals
of Steiner systems, when $a = n-d-1$, cyclic polytopes, when $a=1$ and $d=2c$,
and bi-Cohen-Macaulay simplicial complexes, when $a=0$.

This paper started out partly from an observation that the Tate resolution
(see Section 1) of the famous Horrocks-Mumford bundle on $\PP^{4}$ contains
a $2 \times 5$ matrix of exterior (quadratic) monomials. This paper 
may be considered as studying (complexes of) coherent sheaves on
$\PP^{n-1}$ whose Tate resolution involves a $1 \times N$ matrix of exterior 
monomials (or equivalently a resolution of a monomial ideal in the exterior
algebra). In our investigations we have repeatedly had the benefit of 
computing resolutions over the exterior algebra using Macaulay2 \cite{Mac2}, 
and we express our appreciation of this program.

\section{The BGG correspondence}

We start by recalling some facts about the BGG correspondence 
originating from \cite{BGG}. Our main reference is \cite{EFS}.

\subsection*{Tate resolutions}

Let $V$ be a finite dimensional vector space of dimension $n$ over 
a field $k$. Let $E(V) = \oplus \wedge^i V$ be the exterior algebra and 
for short denote it by $E$. Given a
graded (left) $E$-module $M = \oplus M_i$ we can take a minimal 
projective resolution of $M$
\[ P \, : \quad  \cdots \pil P^{-2} \pil P^{-1} \pil M  \]
where 
\[ P^{-p} = \oplus_{a \in \hele} E(a) \te_k \tilde{V}^{-p}_{-a} . \]

Now the canonical module $\ome$, which is $\Hom_k(E,k)$,
is the injective envelope of $k$. Hence we can take
a minimal injective resolution 
\[ I \, : \,   M \pil I^0 \pil I^1 \pil \cdots  \]
where 
\[ I^p = \oplus_{a \in \hele} \ome(a) \te_k V^p_{-a}. \]

(For $-p < 0$ we put $V^{-p}_{-a}= \tilde{V}^{-p}_{-a-n}$.)
By fixing an isomorphism $k \pil \wedge^n V^*$ where $V^*$ is the dual space
of $V$, we get an isomorphism of $E$ and $\ome(-n)$ as left $E$-modules,
where we have given $V$ degree $1$ and $V^*$ degree $-1$.

We can then join together $P$ and $I$ into an unbounded acyclic complex 
$T(M)$, called the Tate resolution of $M$ 
\[ \cdots \pil \oplus_a \ome(a) \te_k V^p_{-a} \mto{d_p} 
          \oplus_a \ome(a) \te_k V^{p+1}_{-a} \pil \cdots  \]
such that $M$ is $\ker d^0$ and also $\im d^{-1}$.
(One should use $\ome$ instead of $E$ in this complex since 
$\ome$ is the natural thing to use in the framework of Koszul duality
and hence in the BGG correspondence, see \cite{EFS}.)

\subsection*{BGG correspondence}

The terms $T^i$ have natural algebraic geometric interpretations
via the BGG correspondence. 
Let $V$ have a basis $\{e_a\}$ and let $W = V^*$ be the dual space of $V$ 
with dual basis $\{ x_a \}$.
Let $S = S(W)$ be the symmetric algebra on $W$. To $M$ we then associate
a complex of free $S$-modules
\[ L(M) \, : \, \cdots \pil S(i) \te_k M_{i} \mto{\delta^i} S(i+1) \te_k M_{i+1}
  \pil \cdots \]
where 
\[ \delta^i(s \te m) = \sum_a s x_a \te e_a m. \]
If we sheafify $L(M)$ we get a complex of coherent sheaves on the 
projective space $\psp{W}$
\[ \tL(M) \, : \, \cdots \pil \opw(i) \te_k M_i \pil \opw(i+1) \te_k M_{i+1}
   \pil \cdots  . \] 
This, in short, is the BGG correspondence 
between finitely generated graded (left) 
$E$-modules and complexes of coherent sheaves on $\psp{W}$.

\medskip
Suppose $\tL(M)$ has only one non-vanishing cohomology group; a 
coherent sheaf $\gF$. Then the terms of the Tate resolution $T(M)$ 
give the cohomology groups $H^i(\psp{W},\gF(a))$ of $\gF$ 
(for short $H^i\gF(a)$).
More precisely, if $\gF$ is $H^0 \tL(M)$ then 
\begin{equation} T^p(M) = \oplus_{i} \ome(p-i) \te_k H^i \gF(p-i). 
\label{1Tpgen} \end{equation}

Since for a coherent sheaf $\gF$ the cohomology $H^i \gF(a)$ vanishes for
$a \gg 0$ when $i> 0$, we see that for large $p$ 
\begin{equation} T^p(M) = \ome(p) \te_k H^0 \gF(p). 
\label{1Tpcoh} \end{equation}

Conversely, if $M$ is such that $T^p(M)$ is equal to $\ome(p) \te_k
V^p_{-p}$ for large $p$, then the only non-zero cohomology of $\tL(M)$
is in degree $0$ and so $\tL(M)$ gives us a coherent sheaf $\gF$.

\medskip
In general all the $\tL(\ker d^p)[-p]$ for $p$ in $\hele$ 
have the same cohomology, where $[-p]$ denotes the complex shifted 
$p$ steps to the left. Hence $\gF$ is equal to $H^{-p}(\tL(\ker d^p))$
for all $p$. Therefore if we find that the only non-vanishing
cohomology group of $\tL(M)$ is $\gF$ in degree $-p$, we shall think
of $M$ as $\ker d^p$ in  $T$. Then (\ref{1Tpgen}) and (\ref{1Tpcoh})
still hold.

\rem
The BGG correspondence induces an equivalence of 
triangulated categories between the
stable module category of finitely generated graded modules over $E$ and
the bounded derived category of coherent sheaves on $\psp{W}$
\[ \modstab{E} \overset{\tL}{\simeq} D^b(\coh / \psp{W}). \]
\remfin

Due to this remark we may also start with a coherent sheaf $\gF$, and
there will be a module $M$ over $E$ such that $\tL(M)$ only has
non-zero cohomology in degree $0$, equal to $\gF$. Forming the Tate
resolution $T(M)$ we also denote it by $T(\gF)$ and say it is the Tate
resolution of $\gF$.

\subsection*{Duals} Consider $\wedge^n W$ as a module situated in degree
$-n$ and let $M^{\vee}$ be $\Hom_k(M, \wedge^n W)$. Since $\tL(\wedge^n W)$
naturally identifies with the canonical sheaf $\omega_{\psp{W}}$ on 
$\psp{W}$ shifted $n$ places to the left, we see that
\[ \tL(M^{\vee}) = \Hom_k(\tL(M), \omega_{\psp{W}})[n].\]
Hence if $\tL(M)$ has only one nonvanishing cohomology group $\gF$ in 
cohomological degree $p$, then
\begin{equation}
\gext^i(\gF, \omega_{\psp{W}}) = H^{i-p-n}\tL(M^\vee). 
\label{1dual}
\end{equation}

Since $\omega_E$ naturally identifies with $\Hom_k(\omega_E, \wedge^n W)$
we also get that the Tate resolution of $M^{\vee}$ is the dual
$\Hom_k(T(M), \wedge^n W)$ of the Tate resolution of $M$.

\subsection*{Projections} Given a subspace $U \sus W$ we get a projection
$\pi : \psp{W} \projpil \psp{U}$.
If the support of the coherent sheaf $\gF$ does not intersect the
center of projection $\psp{W/U} \sus \psp{W}$ we get a coherent sheaf
$\pi_*(\gF)$ on $\psp{U}$. How is the Tate resolution of $\pi_*\gF$
related to that of $\gF$?
Via the epimorphism $E \pil E(U^*)$ the latter becomes an $E$-module.
It then turns out that the Tate resolution 
\[ T(\pi_* \gF) = \Hom_E(E(U^*),T(\gF)). \]
Note that 
\[ \Hom_E(E(U^*), \ome) = \omega_{E(U^*)}. \]
Hence
\[ T(\pi_*\gF) \, : \, \cdots \pil \oplus_{i} \omega_{E(U^*)}(p-i) \te_k
   H^i \gF(p-i) \pil \cdots . \] 
In particular we see that the cohomology groups $H^i \pi_* \gF(p-i)$ and
$H^i \gF(p-i)$ are equal.

\subsection*{Linear subspaces} If $U \pil W$ is a surjection, we get
an inclusion of linear subspaces $i : \psp{W} \hookrightarrow \psp{U}$. Then by
\cite[1.4 (21)]{Fl} the Tate resolution of $i_*\gF$ is 
\[ \Hom_E(E(U^*), T(\gF)). \]

\section{Simplicial complexes giving coherent sheaves.}
  
\subsection*{The BGG-correspondence applied to simplicial complexes}
Let $\Delta$ be a simplicial complex on the set $[n] = \{1, \ldots, n\}$.
Then we get a monomial ideal $I_\Delta$ in $E$ which is generated by the
monomials $e_{i_1}\cdots e_{i_r}$ such that $\{i_1, \ldots, i_r\}$ is 
not in $\Delta$. Dualizing the inclusion $I_\Delta \sus E(V)$ we get an
exact sequence
\begin{equation}
 0 \pil C_\Delta \pil E(W) \pil (I_\Delta)^* \pil 0. \label{2basis}
\end{equation}
Note that $E(W)$ is a coalgebra and that $C_\Delta$ is the subcoalgebra
generated by all $x_{i_1}\cdots x_{i_r}$ such that $\{i_1, \ldots i_r \}$
is in $\Delta$.

Now think of $\ome = E(W)$ as a left $E(V)$-module; then $C_\Delta$ is
a submodule of $\ome$. Then we can use the BGG correspondence. A natural
question to ask is: When does $\tL(C_\Delta)$ have only one non-vanishing
cohomology group, a coherent sheaf $\gF$? It turns out that this happens
exactly when the Alexander dual simplicial complex $\Delta^*$ is 
Cohen-Macaulay. Let us recall this and some other notions.
\medskip

A simplicial complex $\Delta$ is {\it Cohen-Macaulay} if its 
Stanley-Reisner ring $k[\Delta]$ is a Cohen-Macaulay ring. For more
on this see Stanley's book \cite{St}.

The {\it Alexander dual} $\Delta^*$ of $\Delta$ consists of subsets $F$ of 
$[n]$ such that $[n] - F$ is not a face of $\Delta$. Via the isomorphism
$\ome \iso E(n)$, the submodule $C_{\Delta^*}$ corresponds to the ideal 
$I_{\Delta}$ in $E$. So we get from (\ref{2basis}) an exact sequence
\begin{equation}
0 \pil C_\Delta \pil \omega_E \pil (C_{\Delta^*})^{\vee} \pil 0. 
\label{2basis1}
\end{equation}
Dualizing this we get
\begin{equation}
0 \pil C_{\Delta^*} \pil \omega_E \pil (C_{\Delta})^{\vee} \pil 0. 
\label{2basis2}
\end{equation}

\subsection*{Main theorem}
A coherent sheaf $\gF$ on a projective space is
{\it locally Cohen-Macaulay} of pure dimension $n$
if for all the localizations $\gF_P$ we have depth$\,\gF_P = \dim \gF_P = n$.
This is equivalent to all intermediate cohomology groups $H^i \gF(p)$ vanishing
for $0 < i < n$ when $p$ is large positive or negative. 
It is also equivalent to $\gF$ projecting down to a vector
bundle on ${\bf P}^n$.

Let $c$ be the largest integer such that all $(c-1)$-simplexes of $[n]$ 
are contained in $\Delta$.

\begin{theorem}
\begin{itemize}
\item[a)]
 The complex $\tL(C_\Delta)$ has at most one non-vanishing cohomology group,
a coherent sheaf
$\gF$, if and only if $\Delta^*$ is Cohen-Macaulay. In this case
$\gF$ is $H^{-c} \tL(C_\Delta)$.

\item[b)] $\gF$ is locally Cohen-Macaulay of pure dimension if and only 
if both $\Delta$ and $\Delta^*$ are Cohen-Macaulay.

\item[c)] The support of $\gF$ is contained in a hyperplane if and only if 
$\Delta$ (or equivalently $\Delta^*$) is a cone.

\end{itemize}  \label{2knipp}
\end{theorem}

\begin{proof}
By \cite{ER98} $\dlst$ is Cohen-Macaulay if and only if the associated
ideal of $\dl$ in the symmetric algebra has a linear resolution. 
By \cite[Cor.2.2.2]{AAH99} this happens exactly when $I_\dl$ has a linear
resolution over the exterior algebra. Now note that since $I_\dl$ in $E$
is generated by exterior monomials, in any case a resolution will
have terms
\[ I_\dl \rpil \oplus_{a \geq c+1} E(-a) \te_k \tV_a^1 \rpil 
   \oplus_{a \geq c+2} E(-a) \te_k \tV_a^2 \rpil \cdots \]
with all $\tV^{i}_{c+i}$ non-zero. But then the injective resolution of 
the vector space dual $(I_\dl)^*$ will have ''pure'' terms 
$\omega_E(a) \te_k \tV_{-a}^{a+a_0}$ for $a \gg 0$, meaning $\tL(C_\dl)$ is
a coherent sheaf, if and only if $I_\dl$ has a linear resolution 
from the very start and this then happens exactly when $\dlst$ is 
Cohen-Macaulay.

The fact that $\gF$ is locally Cohen-Macaulay means that the terms
in the Tate resolution are $\omega_E(a) \te_k \tV_{-a}^{a+a_0}$ for
$a \gg 0$ and similarly for $a \ll 0$.

Now by the dual sequences (\ref{2basis1}) and (\ref{2basis2}), the dual
of the Tate resolution of $C_\dl$ is the Tate
resolution of $C_{\dlst}$. Thus we get 
that the condition just stated for the Tate resolution of $\gF$ must mean that
both $\dl$ and $\dlst$ are Cohen-Macaulay.
\medskip

Suppose now the support of $\gF$ is contained in the hyperplane 
$\psp{U} \inpil \psp{W}$
corresponding to a surjection $W \pil U$, where the kernel is generated by
a form $w$ in $W$ defining the hyperplane. Considering $\gF$ as a sheaf
on $\psp U$ denote it by $\gFp$. Then the Tate resolutions are related
by
\[ T(\gF) = \Hom_{E(U^*)}(E, T(\gFp)). \]
Hence the component of $T(\gFp)$ in degree $c$ is $\omega_{E(U^*)}$.
Let the image in $T(\gFp)^c$ of the differential be $C^\prime$.
Then $C_\dl$ is the image of 
\[ \Hom_{E(U^*)}(E,C^\prime) \inpil \Hom_{E(U^*)}(E,\omega_{E(U^*)}) = 
E(W)\] 
and this is again the sum 
$C^\prime + w C^\prime$. Since $C_\dl$ is homogeneous for the multigrading,
we see that $C^\prime$ must also be, and then also $w$, so $w=x_i$ for some
$i$. Then we see that $\dl$ is a cone over the vertex $i$. Since the
argument is clearly reversible, we get c).
\end{proof}

\defi If $\dlst$ is Cohen-Macaulay we denote the corresponding 
coherent sheaf by $\Sdl$.

When $\dl$ and $\dlst$ are both Cohen-Macaulay we say that $\dl$ is
{\it bi-Cohen-Macaulay}.
\defifin

\rem The complex $L(C_\dl)$ is the cellular complex  we get from
$\dl$ by attaching the monomial $x_i$ to the vertex $i$. See \cite{BS}.
\remfin

\begin{proposition} \label{2hilbCM} When $\dlst$ is CM the complex
\[ \tL((C_{\dlst})^{\vee})[-c-1] \, : \,  
\opw(-c-1)^{f^*_{d^*}} \vpil \cdots \vpil \opw(-n) \]
is a resolution of $\gS(\dl)$.
\end{proposition}

\begin{proof} The exact sequence (\ref{2basis1}) gives an exact sequence of
complexes 
\[ 0 \pil \tL(C_\dl) \pil \tL(\omega_E) \pil \tL((C_{\dlst})^{\vee}) \pil 0\]
from which this follows by the long exact cohomology sequence.
\end{proof}

\subsection*{Numerical invariants}
For a simplicial complex $\dl$ on $n$ vertices, let $f_i$ be the number of 
$i$-dimensional simplices. The $f$-polynomial of $\Delta$ is 
\[ f_\dl(t) = 1 + f_0t + f_1t^2 + \cdots + f_{d-1}t^{d}\]
where $d-1$ is the dimension of $\dl$. 

If we form the cone $\cdel$ of $\dl$ over a new vertex, then the
$f$-polynomial of $\cdel$ is 
\[ f_\cdel(t) = (1+t) f_\dl(t). \]

The $f$-polynomial of the Alexander dual $\dlst$ is related to $f$ by
\[ f^*_i + f_{n-i-2} = \binom{n}{i+1}. \]

Note that
the invariants $c^*$ and $d^*$ of $\dlst$ are related to those of $\dl$ by 
\[ c^* + d + 1 = n, \quad c+ d^* + 1 = n. \]

\begin{proposition} Suppose $\dlst$ is Cohen-Macaulay. The Hilbert series
of $\Sdl$ is given by
\[ \sum_k h^0 (\Sdl(k)) t^k = (-1)^{c+1} + (-1)^c f_\dl(-t)/(1-t)^n. \]

If $f_\dl$ is  $(1+t)^{n-s} f$ where $f(1)$ is non-zero, 
then the support of 
$\Sdl$ has dimension $s-1$. \label{3hilb}
\end{proposition}

\begin{proof} The sheaf $\Sdl$ is the cohomology of the complex
\[ \opw(-d)^{f_d} \pil \cdots \pil \opw(-c)^{f_c} \pil \cdots \opw \]
at the term $\opw(-c)^{f_c}$. Since the Hilbert series of 
$\opw(-a)$ is $t^a/(1-t)^n$ we get the proposition by breaking the complex
into short exact sequences and running sheaf cohomology on twists of these.

The statement about the dimension of $\Sdl$ follows by writing $f$
as a polynomial in $(1+t)$.
\end{proof}

There is also another equivalent set of numerical invariants of $\dl$. They
are related to the $f_i$'s by the following polynomial equation
\begin{equation} t^d + f_0t^{d-1} + \cdots + f_{d-1} = 
  (1+t)^d + h_1(1+t)^{d-1} + \cdots + h_d. \label{3fh}
\end{equation}
When $\dl$ is Cohen-Macaulay 
all the $h_i \geq 0$, \cite[II.3]{St}.

There is no geometric interpretation of the $h_i$'s in terms of the
topological realization of $\Delta$. However the following gives a kind
of geometric interpretation of the $h^*_i$'s for a CM simplicial complex
$\dl^*$ in terms of the sheaf $\gS(\Delta)$.

\begin{proposition} If $\dlst$ is CM then in the Grothendieck group of 
sheaves on $\psp{W}$
\begin{equation}
[\gS(\Delta)(c+1)] = h^*_{d^*}[\opnmi] + h^*_{d^* -1}[\go_{{\bf P}^{n-2}}]
+ \cdots + h^*_0[\go_{{\bf P}^{c}}]
\label{2groth}
\end{equation}

More concretely, $\gS_0 = \gS(\dl)(c+1)$ has rank $h^*_{d^*}$ and $\gS_0$ is 
generated by its sections. Take a general map
\[ \opnmi{}^{h^*_{d^*}} \pil \gS_0\]
and let $\gS_1$ be the projection to ${{\bf P}^{n-2}}$ of its cokernel. 
It has rank $h^*_{d^*-1}$ and is generated by its sections. In this way 
we continue.
\end{proposition}

\begin{proof} By Proposition \ref{2hilbCM} there is a resolution
\[ \gS(\dl) \vpil \opw(-c-1)^{f^*_{d^*}} \vpil \cdots \vpil \opw(-n) \]
so the Hilbert series of $\gS_0 = \gS(\dl)(c+1)$ is 
\[ \sum_{i=0}^{d^*} (-t)^i f^*_{d^*-i} /(1-t)^n = 
\sum_{i=0}^{d^*}  h^*_{d^*-i}/(1-t)^{n-i}.\]
 This gives the statement about
the class in the Grothendieck group and so the rank of $\gS_0$ is $h^*_{d^*}$.
Also note by the Tate resolution of $\gS(\dl)$ that $\gS_0$ is $0$-regular
as a coherent sheaf. Consider now the sequence
\[ \opnmi{}^{h^*_{d^*}} \pil \gS_0 \pil \gT_1 \]
where the first is a general map and $\gT_1$ is the cokernel.
Since $\gS_0$ is $0$-regular, $\gT_1$ will also be. Also the Hilbert series
of $\gT_1$ is 
\begin{equation}
\sum_{i=1}^{h^*_{d^*}} h^*_{d^*-i}/(1-t)^{n-i}.
\label{2hilbT}
\end{equation}
Hence letting $\gS_1$ be the projection of $\gT_1$ by a general projection to 
${{\bf P}^{n-2}}$, then since $\gS_1$ and $\gT_1$ have the same cohomology,
$\gS_1$ is $0$-regular with Hilbert series (\ref{2hilbT}). In this way we
may continue.
\end{proof}

\rem We thus see that with larger and larger $c$ we are situated in 
a smaller and smaller part of the Grothendieck group.
\remfin

\section{Bi-Cohen-Macaulay simplicial complexes}

\subsection*{Numerical invariants}
 The basic types of bi-Cohen-Macaulay simplicial complexes
turn out to be the skeletons of 
simplices of various dimensions. So let 
\begin{equation} 
f_{s,c}(t) = 1 + st + \binom{s}{2} t^2 + \cdots + \binom{s}{c} t^c 
\label{3fsc}
\end{equation}
be the $f$-polynomial of the $(c-1)$-dimensional skeleton of the 
$(s-1)$-simplex.

\begin{proposition} If $\dl$ is bi-CM then
\[ f_\dl(t) = (1+t)^{n-s} f_{s,c}(t) \] for some $s$. 
We then say that $\dl$ is 
of type $(n,c,s)$. \label{3fbiCM}
\end{proposition}

\begin{proof}
By (\ref{1dual}) we have that  $\gext^i(\gS(\dl), \omega_{\psp{W}})$ is 
$H^{i+c-n} \tL((C_\dl)^{\vee})$ and by the sequence (\ref{2basis2}) this
identifies with $H^{i+c+1-n} \tL(C_{\dlst})$.

Thus when $\dl$ is bi-CM and so $\gS(\dl)$ is locally Cohen-Macaulay
of dimension $s-1$, then
\[ \gext^{n-s}(\gS(\dl), \omega_{\psp{W}}) = H^{c+1-s} \tL(C_{\dlst}) \]
and the other $\gext$-sheaves vanish. Thus $c+1-s = -c^*$ and since
$c^* + d + 1 = n$ we get $d = n-s + c$.

Now if for a polynomial $f$ we let $c \geq 1$ be the largest integer 
for which 
\[ f(t) = 1 + a + \binom{a}{2} t^2 + \cdots + \binom{a}{c} t^c + \ldots \]
then it is easily seen that $f(t)$ and $(1+t)f(t)$ have the same invariant $c$.
Applying this to $f_{\dl}(t) = (1+t)^{n-s} f(t)$ we see that for the 
polynomial $f$ the degree must be equal to the invariant $c$ and so 
$f(t) = f_{s,c}(t)$. 
\end{proof}

\rem This can of course also rather easily be proven in other ways. 
For instance using the Stanley-Reisner ring $k[\dl]$. Then $\dl$ is
bi-CM iff $k[\dl]$ is CM and has a linear resolution by \cite{ER98}.
By Ex.4.1.17 of \cite{Br-He} it is a simple matter to check that the
$f$-polynomial has the above form.

It can also be deduced numerically by appealing only to the fact that
the $h$-vectors of $\dl$ and $\dlst$ are both non-negative.

%
%
%
%
%

\rem If $\dl$ is Cohen-Macaulay, the terms of the $h$-vector are all
non-negative.  If $\dl$ is bi-CM of type $(n,c,s)$, the terms
$h_{c+1}=h_{c+2}=\cdots =0$.  So the bi-CM simplicial complexes are in
a way numerically extremal in the class of Cohen-Macaulay complexes.
\remfin

\subsection*{Algebraic geometric description of bi-CM simplicial complexes}
Let 
$\Omega^c_{\PP^{n-1}}$ be the sheaf of $c$-differentials on $\PP^{n-1}$.

\begin{proposition}  \label{3skjelett}
a) Let $\dl$ be the $(c-1)$-skeleton of a simplex on $n$ vertices. Then
$\dl$ is bi-CM with $\Sdl = \Omega^c_{\pnmi}$.


b) When $\dl$ is bi-CM of type $(n,c,s)$ 
then $\pi_* \Sdl = \Omega^c_{\psmi}$
where $\pi$ is a projection $\pnmi \dotpil \psmi$ whose center is disjoint from
the support of $\Sdl$.
\end{proposition}

\begin{proof} 
a) When $\dl$ is the $(c-1)$-skeleton of a simplex on $n$ elements
then $\tL(C_\dl)$ is the truncated Koszul complex
\[ \opnmi(-c) \te_k \wedge^c W \pil \cdots \pil \opnmi(-1) 
\te_k W \pil \opnmi. \]
The only cohomology is the kernel of the first map which is 
$\Omega^c_{\psmi}$. Since this is a vector bundle, $\dl$ is bi-CM.

We now prove b). The Tate resolution of $\Sdl$ is
\[ \cdots \pil \omega_E \pil \omega_E(c+1) \te_k H^0 (\Sdl(c+1))
\pil \omega_E(c+2) \te_k H^0 (\Sdl (c+2)) \pil \cdots . \]
For a general subspace $U \sus W$ of dimension $s$
the projection $\pi_* \Sdl$ on $\psmi = {\bf P}(U)$ 
has (minimal) Tate resolution
\[ T(\pi_* \Sdl) = \Hom_{E(V)}(E(U^*), T(\Sdl)). \]

Now $\pi_* \Sdl$ and $\Sdl$ have the same Hilbert series and by 
Proposition \ref{3hilb} this is
the same as the Hilbert series of $\Omega^c_{\psmi}$. When we twist the
latter with $c+1$ its global sections are $\wedge^{c+1} U$.

Now note that if a map 
\[ \omega_{E(U^*)} \pil \omega_{E(U^*)}(c+1) \te_k \wedge^{c+1}U \]
is surjective in degree $-c-1$, 
then it is the map whose graded dual is the unique natural map
\[ \wedge^{c+1} U^* \te_k E(U^*)(-c-1) \pil E(U^*) \]
given by $\wedge^{c+1} U^* \te_k 1 \pil \wedge^{c+1} U^*$. 

Hence the maps
\begin{eqnarray*} 
\omega_{E(U^*)} & \pil & \omega_{E(U^*)}(c+1) \te_k H^0 \pi_* \Sdl(c+1). \\
\omega_{E(U^*)} & \pil & \omega_{E(U^*)}(c+1) \te_k H^0 
\Omega^c_{\PP^{n-1}}(c+1) \end{eqnarray*}
may be identified and so we must have $\pi_* \Sdl = \Omega^c_{\psmi}$.
\end{proof}

\rem In the argument above we actually only used the assumption that
$\dlst$ is CM with $f_\dl$ equal to $(1+t)^{n-s} f_{s,c}(t)$. Only this
thus suffices to conclude that $\dl$ is bi-CM.
\remfin

\subsection*{Topological description of bi-CM simplicial complexes}
The bi-CM simplicial complexes $\dl$ correspond by \cite{ER98} to 
Stanley-Reisner rings $k[\dl]$ which are CM and have a linear 
resolution over the polynomial ring. Since the generators of the
ideal of $k[\dl]$ will have degree $c+1$, 
we say the resolution if $(c+1)$-linear.
 
In \cite{Fr2} R. Fr\"oberg studies Stanley-Reisner rings $k[\dl]$ 
with $2$-linear resolution. When $\dl$ is CM (so $\dl$ is bi-CM with $c=1$)
he shows that $\dl$ is what is called a $(d-1)$-tree. (Strictly speaking
he uses this term only for the $1$-skeleton of $\dl$.) They arise as 
inductively as follows. Start with a $(d-1)$-simplex, then attach $d-1$
simplices, one at a time, by identifying one (and only one) $(d-2)$-face
of $\dl$ with one (and only one) $(d-2)$-face of the simplex to be
attached. This thus describes bi-CM $\dl$ with $c=1$. When $c \geq 2$
things appear to be less tractable as the following example shows.

\eks The following example was first noted in \cite{Re}.
Consider the simplicial complex of dimension $2$ with invariants
$(n,c,s)$ equal to $(6,2,5)$:
 
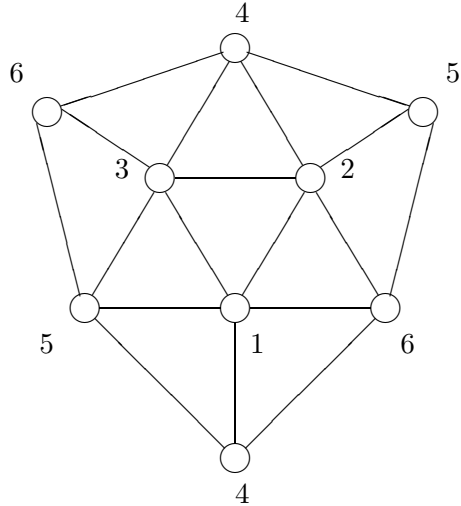
\setlength{\unitlength}{2cm}
\begin{picture}(5,4)
\put(3,0.5){\circle{0.2}}
\put(2,1.5){\circle{0.2}}
\put(4,1.5){\circle{0.2}}
\put(2.5,2.366){\circle{0.2}}
\put(3.5,2.366){\circle{0.2}}
\put(3,3.232){\circle{0.2}}
\put(3,1.5){\circle{0.2}}
\put(1.75,2.8){\circle{0.2}}
\put(4.25,2.8){\circle{0.2}}

\put(2.1,1.5){\line(1,0){0.8}}
\put(3.1,1.5){\line(1,0){0.8}}
\put(2.6,2.366){\line(1,0){0.8}}
\put(2.05,1.583){\line(3,5){0.41}}
\put(3.05,1.583){\line(3,5){0.41}}
\put(2.55,2.453){\line(3,5){0.41}}
\put(3.95,1.583){\line(-3,5){0.41}}
\put(2.95,1.583){\line(-3,5){0.41}}
\put(3.45,2.453){\line(-3,5){0.41}}
\put(3,0.6){\line(0,1){0,8}}
\put(3.07,0.57){\line(1,1){0.86}}
\put(2.93,0.57){\line(-1,1){0.86}}
\put(3.57,2.436){\line(3,2){0.58}}
\put(2.43,2.436){\line(-3,2){0.58}}
\put(4.03,1.58){\line(1,4){0.287}}
\put(1.97,1.58){\line(-1,4){0.287}}
\put(3.08,3.2){\line(3,-1){1.08}}
\put(2.92,3.2){\line(-3,-1){1.08}}

\put(3,0.2){4}
\put(3,3.4){4}

\put(3.1,1.2){1}
\put(1.7,1.2){5}
\put(4.1,1.2){6}
\put(4.4,3){5}
\put(1.5,3){6}
\put(3.7,2.366){2}
\put(2.2,2.366){3}
\end{picture}

This simplicial complex is a triangulation of the real projective
plane.  It is isomorphic to its Alexander dual.  Over any field of
characteristic different from two, it is bi-Cohen-Macaulay.  However,
it has homology in dimension one over $\ZZ/2\ZZ$, so it is not
Cohen-Macaulay over that field.  In particular, it is not shellable.
\eksfin

\section{When are CM-simplicial complexes cones?}

The following proposition gives rise to the problems and results
addressed in this section. In particular we are interested in determining
for which range of invariants $(n,c,s)$  a bi-CM simplicial
complex necessarily is a cone. We give a conjecture for this and prove
the existence of bi-CM simplicial complexes which are not cones in the
whole range of this conjecture.

\begin{proposition}
Let $f$ be a polynomial. Then there exists $e(f)$
such that for $e > e(f)$ if $\dl$ is a CM simplicial complex with
$f_\dl = (1+t)^{e} f$, then $\dl$ is a cone.
\end{proposition}

\begin{proof} The number $c$ of $\dl$ is determined by $f$ and 
the $h^0 (\Sdl(p))$ are also determined by $f$ (Proposition \ref{3hilb}). 
Now by the proof of Proposition \ref{2knipp},
$h^i (\Sdl(c+1-i)) $ is zero for $i > 0$ so $\Sdl$ is $(c+1)$-regular and
is generated by its sections when twisted with $c+1$. Letting $s_1, s_2, 
\ldots, s_a$ be a basis for these sections, there is a surjection
\[ \oplus_{i=1}^a \go_{\pnmi} \te s_i \pil \Sdl(c+1). \]
Now let $b$ be $h^0 (\Sdl(c+2))$.
Then the kernel $K_i$ of each 
\[ H^0 \go_{\pnmi}(1) \te s_i \pil H^0 (\Sdl(c+2)) \]
is at least $(n-b)$-dimensional. If $n > ab$ (which is the case for $e$
sufficiently large), 
the intersection of all the $K_i$ considered as subspaces of
$H^0 \go_{\pnmi}(1)$ is not empty. Thus we get a linear form  $h$ in
$H^0 \go_{\pnmi}(1)$ such that all $h \te s_i$ map to zero. But then
$\Sdl $ is contained in the hyperplane $h=0$ in $\pnmi$ and so $\dl$
is a cone by Proposition \ref{2knipp}.
\end{proof}

We now pose the following.

\begin{problem} For each polynomial $f$ with $f(-1)$ non-zero, 
determine the least number, call it $e(f)$, such that when $\dl$ is
Cohen-Macaulay with $f_\dl = (1+t)^e f$ and {\it not} a cone, then 
$e \leq e(f)$.
\end{problem}

In the case where $f$ is $f_{s,c}$, see (\ref{3fsc}), we propose the following
conjecture for the value of the upper bound of $e = n-s$ when 
$\Delta$ is not a cone.

\begin{conjectureP} Suppose $\dl$ is bi-CM of type
  $(n,c,s)$ and not a cone. Then $n-s  \leq c(s-c-1)$ (or equivalently
  $(c+1)d \leq cn$).
\end{conjectureP}

\begin{conjectureP} Suppose $\gF$ is a non-degenerate
  coherent sheaf on  ${\bf P}^{n-1}$ which projects down to $
  \Omega^c_{\psmi}$ on $\psmi$. Then $n-s \leq c(s-c-1)$
(or equivalently $(c+1)d \leq cn$).
\end{conjectureP}

Clearly Conjecture 2 implies Conjecture 1 by letting $\gF$ be $\Sdl$.
The following shows the existence of non-degenerate coherent sheaves $\gF$
attaining the bound in Conjecture 2 and which cannot be lifted further.

\begin{proposition} The sheaf $\go(-c-1,0)$ on the Segre embedding of
${\bf P}^c \times {\bf P}^{s-c-1}$ in ${\bf P}^{(c+1)(s-c)-1}$ projects 
down to $\Omega^c_{\psmi}$. 
\end{proposition}

Since the Segre embedding is smooth and projectively normal, this line bundle
cannot be lifted further.

\begin{proof} Let us compute the Tate resolution in 
components $c$ and $c+1$. For component $c$ we compute
\[ h^i (\go (-c-1,0)(c-i)) = \left \{ \begin{array}{ll} 0, & i \neq c \\
                                                      1, & i = c 
                                 \end{array}  \right. \]
For component $c+1$ we compute
\[ h^i (\go(-c-1,0)(c+1-i)) = \left \{ \begin{array}{ll} 0, & i > 0 \\
                                    h^0 \go_{{\bf P}^{n-c-1}}(c+1), & i = 0 
                                 \end{array}  \right. \]
Hence components $c$ and $c+1$ of the Tate resolution are
\[ \omega_E \pil \omega_E(c+1) \te_k H^0 \go_{{\bf P}^{n-c-1}}(c+1). \]

 Now note that $h^0 \go_{{\bf P}^{n-c-1}}(c+1)$ is $\binom{n}{c+1}$ 
which again is $h^0 \Omega^c_{\pnmi}(c+1)$. Now the argument proceeds
exactly as in the proof of Proposition \ref{3skjelett} b) and c).
\end{proof}

We shall show the existence of bi-CM simplicial complexes which are
not cones, in the entire range of Conjecture 1.

Let $p$ and $q$ be positive integers. By thinking of a $p \times q$ matrix
we define a {\it vertical} path as a non-decreasing function 
$\alpha : [p] \pil [q]$ and a {\it horizontal} path as a non-decreasing
function $\beta : [q] \pil [p]$. By identifying a path with its graph we may
consider it as a subset of $[p] \times [q]$. We may note that any
horizontal path must intersect any vertical path.

\begin{lemma} Let $F$ be a subset of $[p] \times [q]$. Then either $F$ contains
a horizontal path or the complement $\overline{F}$ contains a vertical path.
\end{lemma}

\begin{proof} We form a partial horizontal path $\beta : [i] \pil [p]$
as follows. Chose $\beta(1)$ minimal such that $(\beta(1), 1)$ is in 
$F$. Then chose $\beta(2) \geq \beta(1)$ minimal such that $(\beta(2), 2)$
is in $F$. Continuing till the process stops gives a path 
$\beta : [i] \pil [p]$. The block $[\beta(i), p] \times [i+1,q] $ can then
contain no element from $F$. Looking at the block 
\[B = [1, \beta(i)-1] \times [1,i], \]
by construction of the path $\beta$, $F \cap B$ does not contain a horizontal
path. By induction $\overline{F} \cap B$ contains a vertical path
$\alpha : [1, \beta(i)-1] \pil [i]$ which can then be completed to a path
all the way down in $\overline{F}$ by picking elements in the block 
$[\beta(i), p] \times [i+1,q] $.
\end{proof}

Let $V$ be the vector space with basis $e_{ij}$ where $i = 1,\ldots, p$ and
$j = 1, \ldots, q$ and fill the matrix with these elements. Let $Y$ be 
the simplicial complex defined by the monomial ideal in $E(V)$ generated
by horizontal path products 
\[ e_{\beta(1), 1} e_{\beta(2), 2} \cdots e_{\beta(q),q} \]
and let $X$ be the simplicial complex defined by the monomial ideal 
generated by the vertical path products.

\begin{proposition}
The facets of $X$ are the complements of the horizontal paths and the facets
of $Y$ are the complements of the vertical paths. In particular $X$ and $Y$
are Alexander dual simplicial complexes.
\end{proposition}

\begin{proof}
The statement about Alexander duals follows from the first statement
because the facets of the Alexander dual of  $Y$ are the subsets of $[p] 
\times [q]$ which are the complements of the indexing sets of the monomial
generators of $Y$, and so this is $X$.

Now the faces of $X$ are precisely the subsets $F$ of $[p] \times [q]$
which do not contain a vertical path. By the previous lemma $\overline{F}$
contains a horizontal path and so $F$ is contained in the complement
of a horizontal path.

If $F$ is the complement of a horizontal path, then $F$ does not contain
a vertical path because any horizontal and vertical paths intersect.
Hence $F$ is a face (in fact a facet) of $X$.
\end{proof}

\begin{lemma}
$X$ and $Y$ are shellable simplicial complexes
\end{lemma}

\begin{proof} The elements of $X$ are complements of horizontal paths
$\beta(1) \cdots \beta(q)$ and we represent them as such. 
We order these lexicographically by letting
$1 \ord 2 \ord \cdots \ord p$. This gives a shelling of $X$.
Let $\alpha $ and $\beta$ be horizontal paths with $\alpha \ord \beta$
and the cardinality of $\alpha \cap \beta $ less than or equal to $q-2$.
If there are at least two values $i$ such that $\alpha(i) < \beta(i)$,
let $l$ be maximal among these and let 
\[ \gamma = \alpha(1) \cdots \alpha(l-1) \beta(l) \cdots \beta(q). \]
If there is only one value $i$ with $\alpha(i) < \beta(i)$ let $i=l$ and
\[ \gamma = \alpha(1) \cdots \alpha(l) \beta(l+1) \cdots \beta(q). \]
Then $\gamma \ord \beta$ and the cardinality of $\alpha \cap \gamma$ is
greater than that of $\alpha \cap \beta$. The argument for $Y$ is similar.
\end{proof}

\begin{theorem} Given $s > c$. For all $n$ in the range 
\[ s \leq n \leq (c+1)(s-c) \]
there exists bi-CM simplicial complexes with invariants $n,c,s$ which are
not cones.
\end{theorem}

There are thus explicit examples in the full range of Conjecture 1.

\begin{proof}
Since $X$ and $Y$ are shellable and Alexander duals,
the Stanley-Reisner ring $k[X]$ is Cohen-Macaulay of dimension 
$pq-q-1$ and has a linear resolution, \cite{ER98}. Considering the
$p \times q$ matrix with entries $x_{ij}$, $k[X]$ is a quotient of 
$k[x_{ij}]$. The $k'th$ {\it normal diagonal} are the positions $(i,j)$
with $i+j = k-1$. We may make some variables on a normal diagonal equal.
Having done so we may form the simplicial complex $X^\prime$ again defined
by the ideal of vertical path products. Then $k[X^\prime]$ is $k[X]$ 
divided out by elements $x_{ij} - x_{i^\prime j^\prime}$ each time we make
$x_{ij}$ equal to $x_{i^\prime j^\prime}$. Making all the elements 
on each normal diagonal equal, call the element on the $k$'th normal diagonal
$x_k$ where $k = 1, \ldots, p+q-1$, we get a simplicial complex $\dl$
on $[p+q-1]$ whose ideal is generated by all monomials $x_{i_1} \cdots x_{i_p}$
where $i_1 < \cdots < i_p$. Thus $\dl$ is the complete $p-2$-dimensional
skeleton of the $p+q-2$-simplex and so $k[\dl]$ is Cohen-Macaulay of 
dimension $p-1$.

  Since we have divided out by $pq-p-q+1$ elements to get from $k[X]$ to 
$k[\dl]$, each time we have cut dimension, and so each element must have
been regular. Therefore $k[X^\prime]$ must be Cohen-Macaulay and have a linear
resolution and so is bi-CM. If we divided out by $m$ elements to get from
$k[X]$ to $k[X^\prime]$ where $m$ is between $0$ and $(p-1)(q-1)$, 
$X^\prime$ has invariants
\[ n = pq -m, \quad d = pq-q-m, \quad c = p-1, \quad s = p+q-1 \]
and so by choosing $p$ and $q$ suitable, we fill up the whole range of the 
theorem.
\end{proof}

Since bi-CM simplicial complexes of type $(n,c,s)$ are Alexander dual to 
bi-CM simplicial complexes of type $(n,s-c-1,s)$ and $\dl$ and $\dlst$
are cones at the same time, we see that if Conjecture 1 is true for  type
$(n,c,s)$ it is true for type $(n,s-c-1,s)$.
The following easy argument shows the conjecture for $c=1$ (and thus
also for $c=s-2$).

\begin{proposition}
If $\dl$ is bi-CM of type $(n,1,s)$ and not a cone, then 
$n-s \leq s-2$.
\end{proposition}

\begin{proof} By Section 3, $\dl$ is constructed as follows. 
Start with a $d$-simplex $F_1$. Attach a $d$-simplex $F_2$ on 
a $(d-1)$-face and continue attaching $F_3, \ldots, F_n$. 
Now $F_1, \ldots, F_s$ are sets of $d$ elements and for each $j$ there is 
$i < j$ such that $F_i \cap F_j$ consists of $d-1$ elements.
But then $\cap_1^s F_i$ contains at least
$d-s+1$ elements and so if $d \geq s$, $\dl$ must be a cone.
Since $d$ is $n-s+1$ this gives the proposition.
\end{proof}

\end{document}